\documentstyle[12pt,amssymb,amscd]{amsart}
\newcommand{\zzz}{{\Bbb Z}}
\newcommand{\jj}{{\frak{g}}}
\newcommand{\hhh}{{\frak{h}}}

\newcommand{\ind}{}
\newenvironment{thm}[1]{\smallskip\ind{\sf #1.}\sl}{\smallskip}
\newcommand{\pr}[1]{\ind{\it #1.\/}}
\newcommand{\Qed}{\qed\smallskip}

\begin{document}
\title{Poincar\'{e} series of the Weyl group of elliptic Lie algebras 
$A^{(1,1)}_1$  and $A^{(1,1) \ast}_1$}
\maketitle
\vspace{-2mm}
\begin{center}
Minoru Wakimoto
\end{center}
\vspace{2mm}
\begin{center}
{\small
Graduate School of Mathematics, Kyushu University, Fukuoka 812-81, Japan}
\end{center}
\vspace{7mm}
\begin{center}
{\bf Abstract}
\end{center}
\begin{center}
\vspace{3mm}
\parbox{12cm}{In this note, we calculate the Poincar\'{e} series of elliptic 
Weyl groups of type $A^{(1,1)}_1$  and $A^{(1,1) \ast}_1$ and show that they 
have a simple expression.  }
\end{center}

\section*{ 0. Introduction}

It is known (see e.g., \cite{B1}, \cite{H1}, \cite{M1}) that the Poincar\'{e} 
series of the Weyl group of a simple Lie algebra has a beautiful expession:
\begin{equation*}
\sum_{w \in W} q^{\ell(w)}  =
\begin{cases}
\prod^{r}_{i=1} \frac{1-q^{m_i +1}}{1-q}, &\qquad \text{for} \,\  \jj ,  \\
\frac{1}{(1-q)^{r}} \prod^{r}_{i=1} \frac{1-q^{m_i +1}}{1-q^{m_i}},
 &\qquad \text{for} \,\  \widehat{\jj} ,  \\
\end{cases}
\end{equation*}
where $r$ is the rank and  $\{m_i; \,\ i= 1, \cdots, r \} $ is the set of 
the exponents, and $\widehat{\jj}$ is the affinization of a 
finite-dimensional simple Lie algebra $\jj$.

It is also known that the notion of extended affine root system was 
introduced by Saito more than a decade ago and that the theory has been 
studied and developed by Saito and other people (e.g., \cite{A1}, \cite{S1}, 
\cite{ST1}, \cite{SL1}, \cite{SL2}, \cite{T1}, \cite{T2}). The Weyl group 
associated to an extended affine root system, although it may be viewed as 
an extension of the affine Weyl group, has a great 
advantage where a Coxeter element comes back alive and plays an important 
role. A Coxeter element is of course an element of finite order (up to a 
certain subgroup in this case), and is closely related to the center $Z$ of 
the Weyl group $W$.

So the interest will naturally arise as to the Poincar\'{e} series of the 
{\it elliptic Weyl group}  $W_F := W/Z$.
This note is a report of an experimental analysis of the Poincar\'{e} series 
of the simplest elliptic Weyl groups of type $A^{(1,1)}_1$ and 
$A^{(1,1)\ast}_1$. Our results may suggest the existence of some kind of 
beautiful formula also in the case of elliptic Weyl groups, although they 
are not Coxeter groups.

The author would thank Dr. K. Iohara for pointing out an 
interesting coincidence, shown in Theorem 3 of \cite{P}, that the 
Weyl group of $A^{(1,1) \ast}_1$ is isomorphic with the Weyl group of 
$SL(2, K)$, where $K$ is a two-dimensional local field.

\section*{1. The Weyl group of the elliptic Lie algebra $A^{(1,1)\ast}_1$. }

The Cartan matrix of $A^{(1,1)\ast}_1$  is
\begin{equation*}
\begin{pmatrix}
  \langle \alpha^{\vee}_0, \,\ \alpha_0 \rangle  
&  \langle \alpha^{\vee}_0, \,\ \alpha_1 \rangle 
&  \langle \alpha^{\vee}_0, \,\ \alpha_{1'} \rangle   \\
  \langle \alpha^{\vee}_1, \,\ \alpha_0 \rangle  
&  \langle \alpha^{\vee}_1, \,\ \alpha_1 \rangle 
&  \langle \alpha^{\vee}_1, \,\ \alpha_{1'} \rangle   \\
  \langle \alpha^{\vee}_{1'}, \,\ \alpha_0 \rangle 
&  \langle \alpha^{\vee}_{1'}, \,\ \alpha_1 \rangle 
&  \langle \alpha^{\vee}_{1'}, \,\ \alpha_{1'} \rangle 
\end{pmatrix} :=
\begin{pmatrix}
       2   &  -1  &  -1   \\
      -4   &   2  &   2   \\
      -4   &   2  &   2
\end{pmatrix},
\end{equation*}
and the Cartan subalgebra $\hhh$ is 5-dimensional
with a basis $\{ d_s, \,\  d_t, \,\ \alpha^{\vee}_i, \,\ (i=0,1,1') \}$, 
where $d_s$ and $d_t$ are elements such that
$\langle d_s, \,\ \alpha_j \rangle = 1$ if $j=1'$, and $=0$ otherwise, and
$\langle d_t, \,\ \alpha_j \rangle = 1$ if $j=1$, and $=0$ otherwise.
Define elements $\Lambda_j$ in the dual space $\hhh^{\ast}$ ($j=0,1,1'$)  by
$\langle \alpha^{\vee}_i, \,\ \Lambda_j \rangle = \delta_{i,j}$ and 
$\langle d_s, \,\ \Lambda_j \rangle = 
\langle d_t, \,\ \Lambda_j \rangle = 0 $, and put 
$\delta := \alpha_0 +2 \alpha_1$ and $\epsilon := \alpha_1 - \alpha_{1'}$.
Then $\{ \Lambda_0, \Lambda_1, \alpha_1, \delta, \epsilon \}$
is a basis of $\hhh^{\ast}$.

Now we consider the Weyl group  $W$ which is, by definition, the subgroup 
of $Aut(\hhh^{\ast})$ generated by reflections $r_i$, for  $i=0,1,1'$, where
\begin{equation*}
r_i(h) := h - \langle h, \,\ \alpha^{\vee}_i \rangle \alpha_i, 
\qquad \text{for all} \,\ 
h \in \hhh^{\ast}.  \tag{1.1}
\end{equation*}

The center $Z$ of $W$ is an infinite cyclic group: \,\ 
$ Z= \langle \sigma^2 \rangle , $
where   $\sigma := r_0 r_1 r_{1'}$   is a Coxeter element.
Actually the action of this element $\sigma^2$ is 
\begin{equation*}
\sigma^2(\Lambda_0) = \Lambda_0 + 4 \epsilon,  \qquad
\sigma^2(\Lambda_1) = \Lambda_1 -  \delta',    \qquad
\sigma^2(\alpha_i) = \alpha_i \,\ \text{for all} \,\ i,   \tag{1.2}
\end{equation*}
where $\delta' := \delta- 2 \epsilon = \alpha_0 + 2 \alpha_{1'}$.
Since $\delta$, $\delta'$ and $\epsilon$ are fixed by all $r_i$, 
(1.2) implies that the action of $\sigma^2$ on $\hhh^{\ast}$ commutes with 
that of all $r_i$, and so $\sigma^2 \in Z$.
One sees also the following:
\begin{equation*}
(r_{0} r_{1} r_{1'})^2
=(r_{1} r_{1'} r_{0})^2
=(r_{1'} r_{0} r_{1})^2
=(r_{0} r_{1'} r_{1})^{-2}
=(r_{1} r_{0} r_{1'})^{-2}
=(r_{1'} r_{1} r_{0})^{-2}.  \tag{1.3}
\end{equation*}

We put $W_F := W/Z$.  By (1.3), one can view  $W_F$ as the group on 
generators $r_i$ , $i=0,1,1'$,  with fundamental relations
\begin{equation*}
 r_i^2  =1,    \tag{1.4a}
\end{equation*}
\begin{equation*}
(r_{i} r_{j} r_{k})^2  =1, \qquad \text{for all permutations} \,\ (i,j,k) \,\ 
\text{of} \,\ (0,1,1').   \tag{1.4b}
\end{equation*}
Note that the condition (1.4b) may be replaced by
\begin{equation*}
r_{i} r_{j} r_{k}=r_{k} r_{j} r_{i} \qquad \text{for all} \,\ (i,j,k).
   \tag{1.4b'}
\end{equation*}

  From (1.4b'), one easily sees that
\begin{equation*}
(r_0 r_1)(r_0 r_{1'}) = (r_0 r_{1'})(r_0 r_1),  \tag{1.5}
\end{equation*}
and that, in view of the exposition of the Weyl group of the affine Lie 
algebra  $A^{(2)}_2$ (e.g., given in Chapter 6 of \cite{K1}), all elements 
in $W_F$ are written uniquely in the form \,\ 
$r_0^i (r_0 r_1)^m  (r_0 r_{1'})^n$, 
where $i= 0,1$, and $ m,n \in \zzz$.

It is also easy to see the following formulas from (1.4) and (1.5):
\begin{equation*}
(r_1 r_{1'})^n = (r_0 r_1)^{-n} (r_0 r_{1'})^n,   \tag{1.6}
\end{equation*}
or equivalently
\begin{equation*}
(r_{1'} r_{1})^n = (r_0 r_1)^{n} (r_0 r_{1'})^{-n}.   \tag{1.6'}
\end{equation*}
for all $n \in \zzz$.

We now consider the length $\ell(w)$ of a minimal expression of an element 
$w \in W_F$ :

\begin{thm}{Lemma 1.1}
\begin{enumerate}
\renewcommand{\labelenumi}{\arabic{enumi})}
\item \quad  For  $w =(r_0 r_1)^m  (r_0 r_{1'})^n$ ( $ m,n \in \zzz$),
\begin{equation*}
\ell(w) = 
\begin{cases}
 2(\vert m \vert + \vert n \vert), & \qquad \text{if} \,\ m,n \ge 0 \,\ 
       \text{or} \,\ m, n \le 0, \\
 2 \text{max} \{\vert m \vert , \vert n \vert \}, & \qquad \text{if} \,\ 
   mn < 0.
\end{cases}
\end{equation*}
\item \quad  For  $w = r_0 (r_0 r_1)^m  (r_0 r_{1'})^n$ ( $ m,n \in \zzz$),
\begin{equation*}
\ell(w) = 
\begin{cases}
 2(m + n)-1, & \qquad \text{if} \,\ m,n \ge 0 \,\ \text{and} \,\ (m,n) 
\ne (0,0), \\
 2(\vert m \vert + \vert n \vert)+1,  & \qquad \text{if} \,\ m,n \le 0,  \\
 2 \text{max} \{\vert m \vert , \vert n \vert \}-1, & \qquad \text{if} \,\ 
  mn < 0     \,\ \text{and} \,\ m+n > 0, \\
 2 \text{max} \{\vert m \vert , \vert n \vert \}+1, & \qquad \text{if} \,\ 
  mn < 0    \,\ \text{and} \,\ m+n \le 0. \\
\end{cases}
\end{equation*}
\end{enumerate}
\end{thm}

\pr{Proof}
1) The cases  $m, n \ge 0$  and $m, n \le 0$ are clear. So we consider the 
case when  $mn < 0$. If $m>0$, $n<0$ and $\vert m \vert \ge \vert n \vert$, 
then
\begin{align*}
w &=(r_0 r_1)^m  (r_0 r_{1'})^n  =(r_0 r_1)^m  (r_0 r_{1'})^{-\vert n \vert}  
  =(r_0 r_1)^{m- \vert n \vert}  \cdot
   (r_0 r_1)^{\vert n \vert}  (r_0 r_{1'})^{-\vert n \vert}       \\
  &=(r_0 r_1)^{m- \vert n \vert} (r_{1'} r_1)^{\vert n \vert}  \qquad \qquad
     \text{by (1.6')},
\end{align*}
so \qquad \qquad \qquad  
$\ell(w) = 2(m- \vert n \vert)+2 \vert n \vert = 2m.$   \\
If $m>0$, $n<0$ and $\vert m \vert \le \vert n \vert$, then
\begin{align*}
w &=(r_0 r_1)^m  (r_0 r_{1'})^n  =(r_0 r_1)^m  (r_0 r_{1'})^{-\vert n \vert}
  =(r_0 r_1)^{m} (r_0 r_{1'})^{-m} \cdot (r_0 r_{1'})^{m-\vert n \vert}   \\
  &=(r_{1'} r_1)^{m} \cdot (r_{1'} r_1)^{\vert n \vert -m}  \qquad \qquad
     \text{by (1.6')},
\end{align*}
so \qquad \qquad \qquad 
$\ell(w) = 2m+ 2(\vert n \vert -m)  = 2 \vert n \vert.$   

Thus we have proved $\ell(w) = \text{max} \{\vert m \vert, \vert n \vert \}$
when $m>0$ and $n<0$.  The case when $m<0$ and $n>0$ is similar, so
1) is proved.     \\
2)  When $m \ge 0$ and $n \ge 0$ and $(m,n) \ne (0,0)$, then
\begin{equation*}
 w =r_0 (r_0 r_1)^m  (r_0 r_{1'})^n  =
\begin{cases}
r_1 (r_0 r_1)^{m-1} (r_0 r_{1'})^n, & \quad  \text{if} \,\ m>0,  \\
r_0 (r_0 r_{1'})^{n} (r_0 r_{1})^m =r_{1'} (r_0 r_{1'})^{n-1} (r_0 r_{1})^m, 
& \quad  \text{if} \,\ n>0, 
\end{cases}
\end{equation*}
so \qquad \qquad \qquad  $\ell(w) = 2(m+n)-1$.    \\
The case when  $m \le 0$ and $n \le 0$   is similar.  

Next we consider the case when $m>0$ and $n<0$.  
If $ m > \vert n \vert$, then
{\allowdisplaybreaks %
\begin{align*}
w &==r_0 (r_0 r_1)^m  (r_0 r_{1'})^{-\vert n \vert}  
  =r_0 (r_0 r_1)^{m- \vert n \vert}  \cdot
   (r_0 r_1)^{\vert n \vert}  (r_0 r_{1'})^{-\vert n \vert}  \\
  &=r_0 (r_0 r_1)^{m- \vert n \vert} (r_{1'} r_1)^{\vert n \vert}  \qquad 
     \qquad    \text{by (1.6')}   \\
  &=r_1 (r_0 r_1)^{m- \vert n \vert -1} (r_{1'} r_1)^{\vert n \vert},
\end{align*}
}
so \qquad \qquad \qquad  
$\ell(w) = 2(m- \vert n \vert)-1+2 \vert n \vert = 2m-1.$   \\
If $ m \le \vert n \vert$, then
\begin{align*}
w &=r_0 (r_0 r_1)^m  (r_0 r_{1'})^{-\vert n \vert}  
 =r_0 \cdot (r_0 r_1)^{m} (r_0 r_{1'})^{-m} 
  \cdot (r_0 r_{1'})^{m-\vert n \vert} \\
 &=r_0 (r_{1'} r_1)^{m} \cdot (r_{1'} r_0)^{\vert n \vert -m}  \qquad \qquad
     \text{by (1.6')},
\end{align*}
so \qquad \qquad \qquad  
$\ell(w) = 2m+1+ 2(\vert n \vert -m)= 2 \vert n \vert +1.$   \\
The case when $m<0$ and $n>0$ is similar.  Thus the Lemma is proved.
\Qed

\begin{thm}{Theorem 1}   For $A^{(1,1)\ast}_1$,   
\begin{equation*}
\sum_{w \in W_F} q^{\ell(w)} = \frac{1- q^3}{(1-q)^3},   \tag{1}
\end{equation*}
\begin{equation*}
\sharp \{ w \in W_F ; \,\ \ell(w) = n \} = 3n   \qquad
\text{if} \,\ n \,\ \text{is a positive integer}.    \tag{2}
\end{equation*}
\end{thm}

\pr{Proof}
By Lemma 1.1,  \,\ $ \sum_{w \in W_F} q^{\ell(w)}$ \,\ is the sum of 
the following six terms:
{\allowdisplaybreaks %
\begin{align*}
A1 &:= \sum_{m,n \ge 0 \,\ \text{or} \,\ m,n \le 0}
q^{2(\vert m \vert + \vert n \vert)},  \qquad \qquad \qquad
A2 := \sum_{mn<0} q^{2 \text{max} \{ \vert m \vert , \vert n \vert \} }, \\
B1 &:= \sum_{m,n \ge 0 \,\ \text{such that} \,\ m+n>0} q^{2(m+n)-1}, 
\qquad \qquad
B2 := \sum_{m,n \le 0} q^{2(\vert m \vert + \vert n \vert)+1},  \\
B3 &:= \sum_{mn<0 \,\ \text{such that} \,\ m+n>0} 
q^{2 \text{max} \{ \vert m \vert , \vert n \vert \}-1 },  \\
B4 &:= \sum_{mn<0 \,\ \text{such that} \,\ m+n \le 0} 
q^{2 \text{max} \{ \vert m \vert , \vert n \vert \}+1 },  
\end{align*}
}
and each one of them is calculated as follows:
{\allowdisplaybreaks %
\begin{align*}
A1 &= 2 \sum_{m,n \ge 0} q^{2(m+n)} -1
= \frac{2}{(1-q^2)^2}-1 = \frac{1+2q^2-q^4}{(1-q^2)^2}, \\
A2 &= 2 \sum_{m>0>n} q^{2 \text{max} \{ m, \vert n \vert \} }  
= 2 \sum_{m,n>0} q^{2 \text{max} \{ m, n \} }    \\
&= 2 \left\{ \sum_{m>0} q^{2m} + 2 \sum_{m>n>0} q^{2m} \right\}  
= 2 \left\{ \frac{q^2}{1-q^2} + \frac{2q^4}{(1-q^2)^2} \right\}, \\
B1 &= \sum_{m,n \ge 0} q^{2(m+n)-1} -q^{-1}
= \frac{1}{q} \left\{ \frac{1}{(1-q^2)^2}-1 \right\}
= \frac{2q-q^3}{(1-q^2)^2}, \\
B2 &=  \sum_{m,n \ge 0} q^{2(m+n)+1} = \frac{q}{(1-q^2)^2},  \\
B3 &= 2 \sum_{m>0>n \,\ \text{such that} \,\ m+n>0}
q^{2 \text{max} \{ m, \vert n \vert \}-1 }  
= 2 \sum_{m>0>n \,\ \text{such that} \,\ m+n>0} q^{2m-1}  \\
&= 2 \sum_{m>k>0} q^{2m-1 } \qquad \qquad \qquad \text{by putting} \,\ 
   n=-k  \\
&= 2 \sum_{m \ge 2} (m-1)q^{2m-1}  
= 2 \cdot \frac{q^3}{(1-q^2)^2},  \\
B4 &= 2 \sum_{m>0>n \,\ \text{such that} \,\ m+n \le 0}
q^{2 \text{max} \{ m, \vert n \vert \}+1 }  
= 2 \sum_{m>0>n \,\ \text{such that} \,\ m+n \le 0} q^{-2n+1}  \\
&= 2 \sum_{k \ge m>0} q^{2k+1 } \qquad \qquad \qquad \text{by putting} \,\ 
n=-k  \\
&= 2 \sum_{k \ge 1} kq^{2k+1}  
= 2 \cdot \frac{q^3}{(1-q^2)^2}.
\end{align*}
}
Then, summing up these six terms, one obtains the formula as desired.
\Qed

\section*{2. The Weyl group of the elliptic Lie algebra $A^{(1,1)}_1$. }

The Cartan matrix of $A^{(1,1)}_1$ is
\begin{equation*}
\begin{pmatrix}
   \langle \alpha^{\vee}_0, \,\ \alpha_0 \rangle
&  \langle \alpha^{\vee}_0, \,\ \alpha_1 \rangle
 & \langle \alpha^{\vee}_0, \,\ \alpha_{0'} \rangle 
&  \langle \alpha^{\vee}_0, \,\ \alpha_{1'} \rangle
  \\
   \langle \alpha^{\vee}_1, \,\ \alpha_0 \rangle
&  \langle \alpha^{\vee}_0, \,\ \alpha_1 \rangle 
 & \langle \alpha^{\vee}_1, \,\ \alpha_{0'} \rangle
&  \langle \alpha^{\vee}_0, \,\ \alpha_{1'} \rangle
  \\
   \langle \alpha^{\vee}_{0'}, \,\ \alpha_0 \rangle
&  \langle \alpha^{\vee}_{0'}, \,\ \alpha_1 \rangle
 & \langle \alpha^{\vee}_{0'}, \,\ \alpha_{0'} \rangle
&  \langle \alpha^{\vee}_{0'}, \,\ \alpha_{1'} \rangle \\
   \langle \alpha^{\vee}_{1'}, \,\ \alpha_0 \rangle
&  \langle \alpha^{\vee}_{1'}, \,\ \alpha_1 \rangle
 & \langle \alpha^{\vee}_{1'}, \,\ \alpha_{0'} \rangle
&  \langle \alpha^{\vee}_{1'}, \,\ \alpha_{1'} \rangle \\
\end{pmatrix} :=
\begin{pmatrix}
       2   &  -2  &   2  &  -2 \\
      -2   &   2  &  -2  &   2 \\
       2   &  -2  &   2  &  -2 \\
      -2   &   2  &  -2  &   2 
\end{pmatrix}.
\end{equation*}
The dual space $\hhh^{\ast}$ of the Cartan subalgebra $\hhh$ is the quotient 
space of the linear span of 
$\{ \Lambda_{\sharp}, \,\ \Lambda_0, \,\ \Lambda_1, \,\ \alpha_i; 
\,\ i=0,0',1,1' \}$   
factored by the equivalence relation 
\begin{equation*}
\alpha_{0}-\alpha_{0'}=\alpha_{1}-\alpha_{1'},    \tag{2.1}
\end{equation*}
where $\Lambda_i$ ($i= \sharp, 0,1$) are defined, so as to be in consistency 
with the equivalence relation (2.1), as follows:
\begin{equation*}
\langle \Lambda_{\sharp}, \,\ \alpha^{\vee}_j \rangle = 1, \qquad 
\text{for all} \,\ j,      \tag{2.2a}
\end{equation*}
\begin{equation*}
\langle \Lambda_i, \,\ \alpha^{\vee}_j \rangle = \delta_{i,j}, \qquad 
\text{for} \,\  i,j=0,1 ,   \tag{2.2b}
\end{equation*}
\begin{equation*}
\langle \Lambda_i, \,\ \alpha^{\vee}_{j'} \rangle = - \delta_{i, 1-j}, 
\qquad \text{for}  \,\ i,j=0,1.   \tag{2.2c}
\end{equation*}
We put 
$\epsilon := \alpha_0 - \alpha_{0'} = \alpha_1 - \alpha_{1'}$,  \,\
$\delta := \alpha_0 + \alpha_1$ and $\delta' := \alpha_{0'} + \alpha_{1'}
= \delta - 2 \epsilon $.

For each $i$,  the simple reflection $r_i \in Aut(\hhh^{\ast})$ is defined by 
the formula (1.1).  Then the action of the element 
$\sigma := r_0 r_{0'} r_1 r_{1'}$  is
\begin{equation*}
\sigma(\Lambda_{\sharp}) = \Lambda_{\sharp} + 2 \epsilon,     \tag{2.3a}
\end{equation*}
\begin{equation*}
\sigma(\Lambda_i) = \Lambda_i - \delta + \epsilon, \qquad \text{for} \,\ 
i=0,1,    \tag{2.3b}
\end{equation*}
\begin{equation*}
\sigma(\alpha_1) = \alpha_1.       \tag{2.3c}
\end{equation*}
So the action of $\sigma$ on $\hhh^{\ast}$ commutes with that of all $r_i$, 
since $\delta$ and $\epsilon$ are fixed by $r_i$.   From this, one sees that 
the Coxeter number of $A^{(1,1)}_1$ is 1, and the 
center $Z$ of the Weyl group $W := \langle r_i ; \,\ i=0,0',1,1'\rangle $ 
is an infinite cyclic group generated by a Coxeter element $\sigma$.
The factor group $W_F := W/Z$  is the group on generators 
$ \{r_i ; \,\ i = 0,0',1,1'\} $   with the fundamental relations
\begin{equation*}
r_i^2 =1,   \tag{2.4a}
\end{equation*}
\begin{equation*}
r_0 r_{0'} r_1 r_{1'} =1.   \tag{2.4b}
\end{equation*}
Note that, from (2.4a) and (2.4b),  one easily has the following:
\begin{equation*}
r_0 r_{0'} = r_{1'} r_{1},   \tag{2.5a}
\end{equation*}
\begin{equation*}
r_0 r_{1'} = r_{0'} r_{1}.   \tag{2.5b}
\end{equation*}

\begin{thm}{Lemma 2.1}
\begin{enumerate}
\renewcommand{\labelenumi}{\arabic{enumi})}
\item \qquad  $(r_0 r_{0'})( r_0 r_{1'}) =(r_0 r_{1'})( r_0 r_{0'})$,
\item \qquad $W_F = \{ (r_0)^i (r_0 r_{0'})^m (r_0 r_{1'})^n; \,\ 
i=0,1 \,\ \text{and} \,\ m, n \in \zzz \}$.
\end{enumerate}
\end{thm}

\pr{Proof}
1) is clear, since 
$$(r_0 r_{0'})(r_0 r_{1'}) =(r_0 r_{0'})( r_{0'} r_{1}) =r_0 r_1$$
 and
$$(r_0 r_{1'})(r_0 r_{0'}) =(r_0 r_{1'})( r_{1'} r_{1}) =r_0 r_1.$$
To prove 2), it suffices to show that 
$r_{0}$, $r_{1}$, $r_{0'}$, and $r_{1'}$  can be written in the form
$ (r_0)^i (r_0 r_{0'})^m (r_0 r_{1'})^n$.
And it is enough to show it only for $r_1$, since this is clear for $r_{0}$,
$r_{0'}$, and $r_{1'}$.  But this is also easy since
\begin{equation*}
r_1 = r_{0'} (r_{0'} r_1)= r_{0'} (r_{0} r_{1'})
= r_{0} (r_{0} r_{0'}) (r_{0} r_{1'}).       \tag{2.6a}
\end{equation*}
\Qed

The following formula, which is an easy consequence from (2.6a),
is also useful in our further calculation:
\begin{equation*}
r_{0}r_{1}= (r_{0} r_{0'}) (r_{0} r_{1'}).       \tag{2.6b}
\end{equation*}

\begin{thm}{Lemma 2.2}
\begin{enumerate}
\renewcommand{\labelenumi}{\arabic{enumi})}
\item \quad  For  $w =(r_0 r_{0'})^m  (r_0 r_{1'})^n$ ( $ m,n \in \zzz$),
$$\ell(w) =  2 \text{max} \{\vert m \vert , \vert n \vert \},  $$
\item \quad  For  $w = r_0 (r_0 r_1)^m  (r_0 r_{1'})^n$ ( $ m,n \in \zzz$),
\begin{equation*}
\ell(w) = 
\begin{cases}
 2 \text{max} \{\vert m \vert , \vert n \vert \}-1, & \qquad \text{if} \,\ 
m+n > 0, \\
 2 \text{max} \{\vert m \vert , \vert n \vert \}+1, & \qquad \text{if} \,\ 
m+n \le 0.
\end{cases}
\end{equation*}
\end{enumerate}
\end{thm}

\pr{Proof}
To prove this Lemma, we make use of the following formulas:
\begin{equation*}
r_{0'}r_{1'}= (r_{0} r_{0'})^{-1} (r_{0} r_{1'}).       \tag{2.7a}
\end{equation*}
\begin{equation*}
r_{1'}r_{0'}= (r_{0} r_{0'}) (r_{0} r_{1'})^{-1}.       \tag{2.7b}
\end{equation*}
In particular from (2.6b), (2.7b) and Lemma 2.1.1), one has
\begin{equation*}
(r_{0}r_{1})^n= (r_{0} r_{0'})^n (r_{0} r_{1'})^n,       \tag{2.6b'}
\end{equation*}
\begin{equation*}
(r_{1'}r_{0'})^n= (r_{0} r_{0'})^n (r_{0} r_{1'})^{-n},       \tag{2.7b'}
\end{equation*}
for all $n \in \zzz$.

1) First consider the case  $m, n \ge 0$. If $m \ge n$, then
\begin{align*}
w &=(r_0 r_{0'})^m (r_0 r_{1'})^n
=(r_0 r_{0'})^{m-n} \cdot (r_0 r_{0'})^n (r_0 r_{1'})^n   \\
&=(r_0 r_{0'})^{m-n} \cdot (r_0 r_{1})^n \qquad \qquad \qquad \text{by 
(2.6b')},
\end{align*}
so \qquad \qquad \qquad $\ell(w) = 2m$.  \\
If $m \le n$, then
\begin{align*}
w &=(r_0 r_{0'})^m (r_0 r_{1'})^n
=(r_0 r_{0'})^{m} (r_0 r_{1'})^m \cdot (r_0 r_{1'})^{n-m}   \\
&=(r_0 r_{1})^{m} \cdot (r_0 r_{1'})^{n-m} \qquad \qquad \qquad \text{by 
(2.6b')},   
\end{align*}
so \qquad \qquad \qquad $\ell(w) = 2n$.  \\
The case $m,n \le 0$ is similar.

Next we consider the case when $m>0$ and $n<0$. For simplicity we put 
$n=-k$.  If $m \ge k$, then
\begin{align*}
w &=(r_0 r_{0'})^m (r_0 r_{1'})^{-k}
=(r_0 r_{0'})^{m-k} \cdot (r_0 r_{0'})^k (r_0 r_{1'})^{-k}   \\
&=(r_0 r_{0'})^{m-k} \cdot (r_{1'} r_{0'})^{k} \qquad \qquad \qquad 
\text{by (2.7b')},
\end{align*}
so \qquad \qquad \qquad $\ell(w) = 2m$.  \\
If $m < k$, then
\begin{align*}
w &=(r_0 r_{0'})^m (r_0 r_{1'})^{-k}
=(r_0 r_{0'})^{m} (r_0 r_{1'})^{-m} \cdot (r_0 r_{1'})^{m-k}   \\
&=(r_{1'} r_{0'})^{m} \cdot (r_{1'} r_{0})^{k-m} \qquad \qquad \qquad 
\text{by (2.7b')},
\end{align*}
so \qquad \qquad \qquad $\ell(w) = 2k$.  \\
The case when $m<0$ and $n>0$ is similar. Thus 1) is proved.  \\
2)  First consider the case $m,n \ge 0$ and $m+n>0$. Then, by (2.6b'), 
one has
\begin{equation*}
w =r_0 (r_0 r_{0'})^m (r_0 r_{1'})^n  =
\begin{cases}
 r_{0} (r_0 r_{0'})^{m-n} (r_0 r_{1})^n 
 =r_{0'} (r_0 r_{0'})^{m-n-1} (r_0 r_{1})^n, &   \text{if} \,\ m > n, \\
 r_{0}  (r_0 r_{1})^n = r_1 (r_0 r_{1})^{n-1}, & \text{if} \,\ m=n,   \\
 r_{0} (r_0 r_{1'})^{n-m} (r_0 r_{1})^{m}
 =r_{1'} (r_0 r_{1'})^{n-m-1} (r_0 r_{1})^{m}, &   \text{if} \,\ m <n,
\end{cases}
\end{equation*}
so \qquad \qquad \qquad  $\ell(w) = 2\text{max}\{m,n \}-1$.    \\
And similarly, if $m,n \le 0$, one has
\qquad  $\ell(w) = 2\text{max}\{ \vert m \vert, \vert n \vert \}+1$.

Next we consider the case when $m>0$ and $n<0$.
For simplicity we put $n=-k$.  If $m \ge k$, then
\begin{align*}
w &=r_0 (r_0 r_{0'})^m (r_0 r_{1'})^{-k}
=r_0 (r_0 r_{0'})^{m-k} \cdot (r_0 r_{0'})^k (r_0 r_{1'})^{-k}   \\
&=r_0 (r_0 r_{0'})^{m-k} \cdot (r_{1'} r_{0'})^{k} \qquad \qquad \qquad 
\text{by (2.7b')},
\end{align*}
so  $\ell(w) = 2m-1$ if $m>k$, and $=2m+1$ if $m=k$.  \\
If $m < k$, then
\begin{align*}
w &=r_0 (r_0 r_{0'})^m (r_0 r_{1'})^{-k}
=r_0 (r_0 r_{0'})^{m} (r_0 r_{1'})^{-m} \cdot (r_0 r_{1'})^{m-k}   \\
&=r_0 (r_{1'} r_{0'})^{m} \cdot (r_{1'} r_{0})^{k-m} \qquad \qquad \qquad 
\text{by (2.7b')},
\end{align*}
so \qquad \qquad \qquad $\ell(w) = 2k+1$.  \\
The case when $m<0$ and $n>0$ is similar. Thus the Lemma is proved. 
\Qed

Now a similar calculation as in the proof of Theorem 1 leads us to the 
following:

\begin{thm}{Theorem 2}   For $A^{(1,1)}_1$,   \\
\begin{equation*}
 \sum_{w \in W_F} q^{\ell(w)} = \frac{(1+q)^2}{(1-q)^2},  \tag{1}
\end{equation*}
\begin{equation*}
\sharp \{ w \in W_F ; \,\ \ell(w) = n \} = 4n   \qquad
\text{if} \,\ n \,\ \text{is a positive integer}.    \tag{2}
\end{equation*}
\end{thm}

\end{document}